\documentclass{amsproc}

\usepackage[swedish,english]{babel}
\usepackage[utf8]{inputenc}
\usepackage[T1]{fontenc}
\usepackage{amsmath,amssymb,amsfonts}
\usepackage{url}
\usepackage[dvipsnames,usenames]{color}
\usepackage{array}
\usepackage[pdftex]{graphicx}
\usepackage{accents}
\usepackage{tabu}

\newcolumntype{L}{>{\displaystyle}l}
\newcolumntype{C}{>{\displaystyle}c}
\newcolumntype{R}{>{\displaystyle}r}

\newcommand{\R}{\mathbb R}
\newcommand{\N}{\mathbb N}

\newcommand{\C}{\mathbb C}

\newcommand{\kilb}{\mathcal K}
\newcommand{\cont}{\mathcal C}
\newcommand{\D}{\mathbb D}
\renewcommand{\Im}{\mathrm{Im}}
\renewcommand{\Re}{\mathrm{Re}}
\def\E{{\mathrm{e}}}
\def\di{\partial} 
\def\til{\widetilde}

\newcommand{\pois}{\mathcal P}
\def\I{\mathfrak{i}}
\newcommand{\diff}{\mathrm{d}}
\renewcommand{\bar}{\overline}

\newcommand{\Dom}{\mathrm{Dom}}

\newcommand{\HN}{Herglotz-Nevan\-linna }

\newcommand{\ntto}{\:\scriptsize{\xrightarrow{\vee\:}}\:}
\newcommand{\cut}{\C\setminus\R}
\newcommand{\cutN}{(\C\setminus\R)^n}

\newcommand{\ie}{\textit{i.e.}\/ } 
\newcommand{\eg}{\textit{e.g.}\/ } 
\newcommand{\cf}{\textit{cf.}\/ } 

\usepackage{bm}
\renewcommand{\vec}{\bm}

\allowdisplaybreaks

\theoremstyle{definition}
\newtheorem{define}{Definition}[section]
\newtheorem{example}[define]{Example}
\newtheorem{remark}[define]{Remark}

\theoremstyle{plain}
\newtheorem{lemma}[define]{Lemma}
\newtheorem{thm}[define]{Theorem}
\newtheorem{prop}[define]{Proposition}
\newtheorem{coro}[define]{Corollary}

\numberwithin{equation}{section}

\begin{document}

\title[Characterizations using positive semi-definite functions]{Characterizations of \HN functions using positive semi-definite functions and the Nevanlinna kernel in several variables}

\author{Mitja Nedic}
\address{Mitja Nedic, Department of Mathematics and Statistics, University of Helsinki, PO Box 68, FI-00014 Helsinki, Finland, orc-id: 0000-0001-7867-5874}
\curraddr{}
\email{mitja.nedic@helsniki.fi}
\thanks{\textit{Key words.} \HN functions, positive semi-definite functions, Poisson-type functions, Nevanlinna kernel. 
}

\subjclass[2010]{32A26, 32A99.}

\date{2018-12-19} 

\begin{abstract}
In this paper, we give several characterizations of \HN functions in terms of a specific type of positive semi-definite functions called Poisson-type functions. This allows us to propose a multidimensional analogue of the classical Nevanlinna kernel and a definition of generalized Nevanlinna functions in several variables. Furthermore, a characterization of the symmetric extension of a \HN function is also given. The subclass of Loewner functions is also discussed, as well as an interpretation of the main result in terms of holomorphic functions on the unit polydisk with non-negative real part.
\end{abstract}

\maketitle

\section{Introduction}
\label{sec:intro}

Let us denote by $\C^{+n}$ the poly-upper half-plane in $\C^n$, \ie
$$\C^{+n}:= (\C^+)^n = \big\{\vec{z}\in\C^n \,\big |\,\forall j=1,2,\ldots, n:  \Im[z_j]>0 \big\}.$$
Here, we consider the following class of functions.
\begin{define}\label{def:hn_functions}
A function $q \colon \C^{+n} \to \C$ is called a \emph{\HN function} if it is holomorphic and has non-negative imaginary part.
\end{define}
This is a well-studied class of functions, appearing \eg in \cite{AglerEtal2012,AglerEtal2016,LugerNedic2017,LugerNedic2019,LugerNedic2020,Nedic2019,Savchuk2006,Vladimirov1969,Vladimirov1979}. Particularly when considering such functions of one variable, they appear in many areas of complex and functional analysis and numerous applications. Some examples of these include the moment problem \cite{Akhiezer1965,Nevanlinna1922,Simon1998}, the theory of Sturm-Liouville problems and perturbations \cite{Aronszajn1957,AronszajnBrown1970,Donoghue1965,KacKrein1974}, when deriving physical bounds for passive systems \cite{BernlandEtal2011} or as approximating functions in certain convex optimization problems \cite{IvanenkoETAL2019a,IvanenkoETAL2020}. Applications concerning functions of several variables appear \eg when considering operator monotone functions \cite{AglerEtal2012} or with representations of multidimensional passive systems \cite{Vladimirov1979}.

\HN functions admit a powerful integral representation theorem, \cf Theorem \ref{thm:intRep_Nvar}, that characterizes this class of functions in terms of a real number $a$, a vector $\vec{b} \in [0,\infty)^n$ and a positive Borel measure $\mu$ on $\R^n$ satisfying two conditions. The one-variable version this representation is considered a classical result attributed to Nevanlinna \cite{Nevanlinna1922} and Cauer \cite{Cauer1932}, see also \cite{KacKrein1974}. The multi-dimensional case of this representation was considered in \eg \cite{LugerNedic2017,LugerNedic2019,Vladimirov1969,Vladimirov1979}.

For \HN functions of one variable, it follows as an immediate consequence of the integral representation theorem mentioned above that a holomorphic function $q\colon\C^+ \to \C$ has non-negative imaginary part if and only if the function
$$(z,w) \mapsto \frac{q(z) - \bar{q(w)}}{z-\bar{w}}$$
is positive semi-definite on $\C^+ \times \C^+$, see also Section \ref{sec:main_thm}. This result gives rise to the so-called Nevanlinna kernel, see also Section \ref{subsec:Nevanlinna_kernel}, which plays a fundamental role in \eg the theory of generalized Nevanlinna functions and their applications \cite{DijksmaShondin2000,KreinLanger1977,KurasovLuger2011,LangerLangerSasvari2004}.

The main goal of this paper is to first give a characterization when a function $q\colon\C^{+n} \to \C$ is a \HN function by investigating the difference $q(\vec{z}) - \bar{q(\vec{w})}$ for $\vec{z},\vec{w} \in \C^{+n}$, \cf Theorem \ref{thm:positive_semi_decomposition}. Here, a crucial role is played by Poisson-type functions which are positive semi-definite functions of a very particular type, see Section \ref{subsec:poisson_type}. This result then allows us to derive an analogous characterization to the one mentioned above using a suitable analogue of the Nevanlinna kernel in several variables, \cf Theorem \ref{thm:Nevan_kernel_Nvar}.

The structure of the paper is as follows. After the introduction in Section \ref{sec:intro} we recall the integral representation theorem for \HN function along with the necessary corollaries in Section \ref{sec:hn_basics}. Required prerequisites concerning positive semi-definite functions are collected in Section \ref{sec:psd}, with Section \ref{subsec:poisson_type} focusing on Poisson-type functions. The main results of the paper are presented in Section \ref{sec:main_thm}, with Section \ref{subsec:Nevanlinna_kernel} focusing on the Nevanlinna kernel in several variables, Section \ref{subsec:symmetric_extension} presenting an extension of the main result to symmetric extensions of \HN functions and Section \ref{subsec:loewner} discussing the subclass of Loewner functions. Finally, Section \ref{sec:polydisk} interprets the main result of the paper in the language of holomorphic functions on the unit polydisk with non-negative real part.

\section{The integral representation theorem}
\label{sec:hn_basics}

\HN functions are primarily studied via their corresponding integral representation theorem, the statement of which requires us to first introduce some notation. Given ambient numbers $z \in \cut$ and $t \in \R$, consider first the expressions
$$\begin{array}{RCL}
N_{-1}(z,t) & := & \frac{1}{2\,\I}\left(\frac{1}{t - z} - \frac{1}{t - \I}\right), \\[0.35cm]
N_{0}(z,t) & := & \frac{1}{2\,\I}\left(\frac{1}{t - \I} - \frac{1}{t_j + \I}\right), \\[0.35cm]
N_{1}(z,t) & := & \frac{1}{2\,\I}\left(\frac{1}{t + \I} - \frac{1}{t - \bar{z}}\right).
\end{array}$$
Note that $N_0(z,t) \in \R$ and is independent of $z \in \cut$, while
$$\bar{N_{-1}(z,t)} = N_{1}(z,t)$$
for all $z \in \cut$ and $t \in \R$. Define now the kernel $K_n \colon \cutN \times \R^n \to \C^n$ as
\begin{multline}\label{eq:kernel_Kn}
    K_n(\vec{z},\vec{t}) := \I\bigg(2\prod_{j=1}^n\big(N_{-1}(z_j,t_j)+N_{0}(\I,t_j)\big)-\prod_{j=1}^nN_{0}(\I,t_j)\bigg) \\
    = \I\left(\frac{2}{(2\,\I)^n}\prod_{\ell=1}^n\left(\frac{1}{t_\ell-z_\ell}-\frac{1}{t_\ell+\I}\right)-\frac{1}{(2\,\I)^n}\prod_{\ell=1}^n\left(\frac{1}{t_\ell-\I}-\frac{1}{t_\ell+\I}\right)\right).
\end{multline}

We may now recall the integral representation theorem for \HN functions of several variables, \cf \cite[Thm. 4.1]{LugerNedic2019}.

\begin{thm}\label{thm:intRep_Nvar}
A function $q\colon \C^{+n} \to \C$ is a \HN function if and only if $q$ can be written as
\begin{equation}\label{eq:intRep_Nvar}
q(\vec{z}) = a + \sum_{\ell=1}^nb_\ell z_\ell + \frac{1}{\pi^n}\int_{\R^n}K_n(\vec{z},\vec{t})\diff\mu(\vec{t}),
\end{equation}
where $a \in \R$, $\vec{b} \in [0,\infty)^n$ and $\mu$ is a positive Borel measure on $\R^n$ satisfying the growth condition
\begin{equation}
\label{eq:measure_growth}
\int_{\R^n}\prod_{\ell=1}^n\frac{1}{1+t_\ell^2}\diff\mu(\vec{t}) < \infty
\end{equation}
and the Nevanlinna condition
\begin{equation}
\label{eq:measure_Nevan}
\sum_{\substack{\vec{\rho} \in \{-1,0,1\}^n \\ -1\in\vec{\rho} \wedge 1\in\vec{\rho}}}\int_{\R^n}\prod_{j=1}^nN_{\rho_j}(z_j,t_j)\diff\mu(\vec{t}) = 0
\end{equation}
for all $\vec{z} \in \C^{+n}$. Furthermore, for a given function $q$, the triple of representing parameters $(a,\vec{b},\mu)$ is unique.
\end{thm}

An important corollary of this results is the description of the growth of a \HN function at infinity in a Stoltz domain. Generally, an \emph{(upper) Stoltz domain} with centre $t_0 \in \R$ and angle $\theta \in (0,\frac{\pi}{2}]$ is the set
$$\{z \in \C^+~|~\theta \leq \arg(z-t_0) \leq \pi-\theta\}.$$
The symbol $z \ntto \infty$ then denotes the limit $|z| \to \infty$ in any Stoltz domain with centre $0$. The growth of a \HN function is then captured by the vector $\vec{b}$, as it holds for any $\ell \in \{1,\ldots,n\}$ that
\begin{equation}
    \label{eq:b_parameters}
    b_\ell  = \lim\limits_{z_\ell \ntto \infty}\frac{q(\vec{z})}{z_\ell},
\end{equation}
see \eg \cite[Cor. 4.6(iv)]{LugerNedic2019}. In particular, the above limit is independent of the entries of the vector $\vec{z}$ at the non-$\ell$-th positions.

\section{Positive semi-definite functions}
\label{sec:psd}

Consider the following definition, \cf \cite[pg. 3002]{AglerEtal2016}.

\begin{define}\label{def:positive_semi_def_function}
Let $\Omega \subseteq \C^n$ be an open set. A function $F\colon \Omega \times \Omega \to \C$ is called \emph{positive semi-definite} if for all $m \geq 1$ and every choice of $m$ vectors $\vec{z}_1,\vec{z}_2,\ldots,\vec{z}_m \in \Omega$ and $m$ numbers $c_1,c_2,\ldots,c_m \in \C$ it holds that
\begin{equation}
    \label{eq:positive_semi_def_function}
    \sum_{i,j = 1}^mF(\vec{z}_i,\vec{z}_j)\,c_i\,\bar{c}_j \geq 0.
\end{equation}
\end{define}

Note that we do not impose any regularity (or analyticity) on the function $F$, although the majority of the positive semi-definite functions that will be considered in this paper will be holomorphic in the first $n$ variables and anti-holomorphic in the second $n$ variables.

Alternative, one can say that a function $F\colon \Omega \times \Omega \to \C$ is positive semi-definite if for all $m \geq 1$ and every choice of $m$ vectors $\vec{z}_1,\vec{z}_2,\ldots,\vec{z}_m \in \Omega$ it holds that the matrix
$$\left[\begin{array}{cccc}
F(\vec{z}_1,\vec{z}_1) & F(\vec{z}_1,\vec{z}_2) & \ldots & F(\vec{z}_1,\vec{z}_m) \\
F(\vec{z}_2,\vec{z}_1) & F(\vec{z}_2,\vec{z}_2) & \ldots & F(\vec{z}_2,\vec{z}_m) \\
\vdots & \vdots & \ddots & \vdots \\
F(\vec{z}_m,\vec{z}_1) & F(\vec{z}_m,\vec{z}_2) & \ldots & F(\vec{z}_m,\vec{z}_m) 
\end{array}\right]_{m \times m}$$
is a positive semi-definite matrix.

\begin{remark}
We note here that what we here call positive semi-definite \emph{functions} are sometimes referred to as positive semi-definite \emph{kernels}, \eg in \cite{Aronszajn1950,BuescuPaixao2014,Sasvari1994}. Our terminology stems from \eg \cite[pg. 3002]{AglerEtal2016} and refers to functions of $2n$ complex variables. These functions should not be confused with positive semi-definite functions in the sense of \eg \cite[Def. 1.3.1]{Sasvari1994}, which refers to functions of one real variable.
\end{remark}

The following elementary properties now follow from the definition and will be used later on, see \eg \cite[Sec. 1.3]{Sasvari1994} for a proof.

\begin{lemma}\label{lem:positive_semi_def_properties}
If $F_1$ and $F_2$ are two positive semi-definite functions on $\Omega \times \Omega$, then the functions $F_1 + F_2$ and $F_1F_2$ are also positive semi-definite. Furthermore, for any positive semi-definite function $F$ on $\Omega \times \Omega$, it holds that $F(\vec{z},\vec{z}) \geq 0$ and
$$F(\vec{z},\vec{w}) = \bar{F(\vec{w},\vec{z})}$$
for any $\vec{z},\vec{w} \in \Omega$.
\end{lemma}

\begin{example}
\label{ex:psd_simple}
Let $\Omega \subseteq \C^n$ and let $f\colon \Omega \to \C$ be a holomorphic function. Then, the function
$$F\colon (\vec{z},\vec{w}) \mapsto f(\vec{z})\bar{f(\vec{w})}$$
is positive semi-definite on $\Omega \times \Omega$. Indeed, if $m \geq 1$,  $\vec{z}_1,\vec{z}_2,\ldots,\vec{z}_m \in \Omega$ and $c_1,c_2,\ldots,c_m \in \C$, then
$$\sum_{i,j = 1}^mF(\vec{z}_i,\vec{z}_j)\,c_i\,\bar{c}_j = \sum_{i,j = 1}^mf(\vec{z}_i)c_i\,\bar{f(\vec{z}_j)}\bar{c}_j = |f(\vec{z}_1)c_1 + \ldots + f(\vec{z}_m)c_m|^2  \geq 0.$$
Note that the last equality above holds by a combinatorial expansion of the square of an absolute value of a sum of complex numbers, \ie if $m \in \N$ and $\zeta_1,\ldots,\zeta_m \in \C$, then
\begin{multline*}
    |\zeta_1 + \ldots + \zeta_m|^2 = (\zeta_1 + \ldots + \zeta_m)(\bar{\zeta}_1 + \ldots + \bar{\zeta}_m) \\ = \zeta_1\bar{\zeta}_1 + \zeta_1\bar{\zeta}_2 + \zeta_2\bar{\zeta}_1 + \ldots + \zeta_m\bar{\zeta}_m = \sum_{i,j=1}^m\zeta_i\bar{\zeta}_j
\end{multline*}
as desired.\hfill $\lozenge$
\end{example}

\subsection{Poisson-type functions}
\label{subsec:poisson_type}

In later sections, we will mostly look at functions on the poly-upper half-plane. There, the following class of functions is of major importance.

\begin{define}\label{def:poisson_type}
A function $F\colon \C^{+n} \times \C^{+n} \to \C$ is called a \emph{Poisson-type function} if there exists a positive Borel measure $\mu$ on $\R^n$ satisfying the growth condition \eqref{eq:measure_growth} such that
$$F(\vec{z},\vec{w}) = \frac{1}{\pi^n}\int_{\R^n}\prod_{\ell = 1}^n\frac{1}{(t_\ell - z_\ell)(t_\ell - \bar{w_\ell})}\diff\mu(\vec{t})$$
for all $\vec{z},\vec{w} \in \C^{+n}$. In this case, we say that the function $F$ is given by the measure $\mu$.
\end{define}

There are several important remarks to be made here. First, the name of these functions refers to their connection to the Poisson kernel of $\C^{+n}$ which will become apparent in the proof of Lemma \ref{lem:Stieltjes_for_Poisson} later on. Second, the above definition makes sense as the assumption that the measure $\mu$ satisfies the growth condition \eqref{eq:measure_growth} assure that the integral
$$\int_{\R^n}\prod_{\ell = 1}^n\frac{1}{(t_\ell - z_\ell)(t_\ell - \bar{w_\ell})}\diff\mu(\vec{t})$$
is well-defined for $\vec{z},\vec{w} \in \C^{+n}$. Third, the normalizing factor $\pi^{-n}$ is present so that the function given by the Lebesgue measure $\lambda_{\R}$ equals $2\,\I(z-\bar{w})^{-1}$. Finally, it is not immediately clear from Definition \ref{def:poisson_type} whether the correspondence between a function $F$ and a measure $\mu$ constitutes a bijection, though this will turn out to be the case later on, \cf Lemma \ref{lem:Stieltjes_for_Poisson}.

We will now show three elementary, but important, properties of Poisson-type functions.

\begin{lemma}
\label{lem:positive_semi_poisson}
Let $F$ be a Poisson-type function given by a measure $\mu$. Then, the function $F$ is positive semi-definite on $\C^{+n} \times \C^{+n}$.
\end{lemma}

\proof
Let $m \in \N$ be arbitrary and choose freely $m$ vectors $\vec{z}_1,\vec{z}_2,\ldots,\vec{z}_m \in \C^{+n}$ and $m$ numbers $c_1,c_2,\ldots,c_m \in \C$. In this case, we calculate that
$$\begin{array}{RCL}
\multicolumn{3}{L}{\sum_{i,j=1}^m F(\vec{z}_i,\vec{z}_j)\,c_i\,\bar{c}_j = \frac{1}{\pi^n}\int_{\R^n}\sum_{i,j=1}^m\bigg[c_i\,\bar{c}_j\prod_{\ell=1}^{n}\frac{1}{(t_\ell - (\vec{z}_i)_\ell)(t_\ell - \bar{(\vec{z}_j)}_\ell)}\bigg]\diff\mu(\vec{t})} \\[0.6cm]
~~~~~~~ & = & \frac{1}{\pi^n}\int_{\R^n}\sum_{i,j=1}^m\bigg[c_i\prod_{\ell=1}^{n}\frac{1}{t_\ell - (\vec{z}_i)_\ell}\,\cdot\,\bar{c_j \prod_{\ell=1}^{n}\frac{1}{t_\ell - (\vec{z}_j)_\ell}}\bigg]\diff\mu(\vec{t}) \\[0.6cm]
~ & = & \frac{1}{\pi^n}\int_{\R^n}\bigg|\sum_{i=1}^m\bigg[c_i\prod_{\ell=1}^{n}\frac{1}{t_\ell - (\vec{z}_i)_\ell}\bigg]\,\bigg|^2\diff\mu(\vec{t}) \geq 0,
\end{array}$$
where the last equality follows by the same combinatorial argument used in Example \ref{ex:psd_simple}. This finishes the proof.
\endproof

\begin{lemma}
\label{lem:poisson_growth}
Let $F$ be a Poisson-type function given by a measure $\mu$. Then, it is holomorphic in its first $n$ variables, anti-holomorphic in its second $n$ variables and satisfies the growth condition that for every $j \in \{1,\ldots,n\}$ we have
\begin{equation}
    \label{eq:positive_semi_growth}
    \lim\limits_{z_j \ntto \infty}F(\vec{z},\vec{w}) = \lim\limits_{w_j \ntto \infty}F(\vec{z},\vec{w}) = 0.
\end{equation}
\end{lemma}

\proof
Fix first an arbitrary $\vec{w} \in \C^{+n}$ and consider the function $\vec{z} \mapsto F(\vec{z},\vec{w})$ on $\C^{+n}$. This function is holomorphic as the kernel
$$\prod_{\ell = 1}^n\frac{1}{(t_\ell - z_\ell)(t_\ell - \bar{w_\ell})}$$
is holomorphic in the $\vec{z}$-variables for every $\vec{t} \in \R^n$ while the function
$$\vec{z} \mapsto \frac{1}{\pi^n}\int_{\R^n}\prod_{\ell = 1}^n\frac{1}{|t_\ell - z_\ell||t_\ell - \bar{w_\ell}|}\diff\mu(\vec{t})$$
is locally uniformly bounded on compact subsets of $\C^{+n}$ due to fact that the measure $\mu$ satisfies the growth condition \eqref{eq:measure_growth}. An analogous argument may now be repeated to show that, for every $\vec{z} \in \C^{+n}$, the function $\vec{w} \mapsto F(\vec{z},\vec{w})$ is anti-holomorphic on $\C^{+n}$.

To prove that the function $F$ satisfies the growth conditions \eqref{eq:positive_semi_growth}, it suffices to only consider the case when $z_1 \ntto \infty$, as all other cases may be treated analogously. In this case, we note, for $\vec{z},\vec{w} \in \C^{+n}$, that
\begin{multline*}
    |F(\vec{z},\vec{w})| \leq \frac{1}{\pi^n}\int_{\R^n}\prod_{j=1}^n\frac{1}{|t_j-z_j||t_j - \bar{w}_j|}\diff\mu(\vec{t}) \\ = \frac{1}{\pi^n}\int_{\R^n}\left|\frac{t_1-\I}{t_1-z_1}\right|\frac{1}{|t_1-\I||t_1-\bar{w}_1|}\prod_{j=2}^n\frac{1}{|t_j-z_j||t_j - \bar{w}_j|}\diff\mu(\vec{t}).
\end{multline*}
As $z_1 \ntto \infty$, we may assume that $\Im[z_1] > 1$, yielding that
$$\left|\frac{t_1-\I}{t_1-z_1}\right| \leq 1$$
for all $t_1 \in \R$. Furthermore, the function
$$\vec{t} \mapsto \frac{1}{|t_1-\I||t_1-\bar{w}_1|}\prod_{j=2}^n\frac{1}{|t_j-z_j||t_j - \bar{w}_j|}$$
is integrable with respect to the measure $\mu$ on $\R^n$ as $\mu$ satisfies the growth condition \eqref{eq:measure_growth}. Thus, by Lebesgue's dominated convergence theorem, we have that
$$\lim\limits_{z_1 \ntto \infty}F(\vec{z},\vec{w}) = \frac{1}{\pi^n}\int_{\R^n}\lim\limits_{z_1 \ntto \infty}\prod_{\ell=1}^{n}\frac{1}{(t_\ell - z_\ell)(t_\ell - \bar{w}_\ell)}\diff\mu(\vec{t}) = 0.$$
Note now that the function $F$ is positive semi-definite by Lemma \ref{lem:positive_semi_poisson}, implying, by Lemma \ref{lem:positive_semi_def_properties}, that $F(\vec{z},\vec{w}) = \bar{F(\vec{w},\vec{z})}$ for any $\vec{z},\vec{w} \in \C^{+n}$. Hence, 
$$\lim\limits_{w_j \ntto \infty}F(\vec{z},\vec{w}) = \lim\limits_{w_j \ntto \infty}\bar{F(\vec{w},\vec{z})} = 0.$$
This finishes the proof.
\endproof

\begin{lemma}
\label{lem:Stieltjes_for_Poisson}
Let $F$ be a Poisson-type function given by a measure $\mu$. Then, the measure $\mu$ may be reconstructed form the function $F$ via the Stieltjes inversion formula. More precisely, let $\psi\colon \R^n \to \R$ be a $\cont^1$-function for which there exists some constant $C \geq 0$ such that $|\psi(\vec{x})| \leq C\prod_{j=1}^n(1+x_j^2)^{-1}$ for all $\vec{x} \in \R^n$. Then, it holds that
\begin{equation}
    \label{eq:Stieltjes_for_Poisson}
    \lim\limits_{\vec{y} \to \vec{0}^+}\int_{\R^n}\psi(\vec{x})\cdot\prod_{j=1}^ny_j\cdot F(\vec{x}+\I\,\vec{y},\vec{x}+\I\,\vec{y})\diff\vec{x} = \int_{\R^n}\psi(\vec{t})\diff\mu(\vec{t}).
\end{equation}
\end{lemma}

\proof
We only need to observe that for a Poisson-type function $F$, it holds that
$$\prod_{j=1}^ny_j\cdot F(\vec{x}+\I\,\vec{y},\vec{x}+\I\,\vec{y}) = \frac{1}{\pi^n}\int_{\R^n}\pois_n(\vec{x}+\I\,\vec{y},\vec{t})\diff\mu(\vec{t}),$$
where $\pois_n$ denotes the Poisson kernel of the poly-upper half-plane which, we recall, is defined for $\vec{z} \in \C^{+n}$ and $\vec{t} \in \R^n$ as
$$\pois_n(\vec{z},\vec{t}) := \prod_{j=1}^n\frac{\Im[z_j]}{|t_j-z_j|^2}.$$
The statement then follows immediately from the properties of the Poisson kernel as in the proof of the Stieltjes inversion formula for \HN functions of several variables, see \cite[pg. 1197]{LugerNedic2019}.
\endproof

Recall now that the representing measure of a \HN function satisfies the Nevanlinna condition \eqref{eq:measure_Nevan} in addition to the growth condition \eqref{eq:measure_growth}. As satisfying condition \eqref{eq:measure_growth} suffices for a positive Borel measure on $\R^n$ to define a Poisson-type function, it is apparent that Poisson-type function given by representing measures of \HN functions constitute a smaller subclass. This is reflected by the following condition.

\begin{lemma}
\label{lem:poisson_pluriharmonic}
Let $F$ be a Poisson type function given by a measure $\mu$. Then, the function
\begin{equation}
    \label{eq:poisson_pluriharmonic}
    \vec{z} \mapsto \prod_{j=1}^n\Im[z_j]\cdot F(\vec{z},\vec{z})
\end{equation}
is pluriharmonic on $\C^{+n}$ if and only if the measure $\mu$ satisfies the Nevanlinna condition \eqref{eq:measure_Nevan}.
\end{lemma}

\proof
We have already noted in the proof of Lemma \ref{lem:Stieltjes_for_Poisson} that
$$\prod_{j=1}^n\Im[z_j]\cdot F(\vec{z},\vec{z}) = \frac{1}{\pi^n}\int_{\R^n}\pois_n(\vec{z},\vec{t})\diff\mu(\vec{t}).$$
The statement of the lemma now follows from \cite[Prop. 5.2 and 5.3]{LugerNedic2019}.
\endproof

\begin{example}
Consider the Poisson-type function $F$ given by the measure $\pi^2\delta_{(0,0)}$, where $\delta_{(0,0)}$ denotes the Dirac measure at $(0,0) \in \R^2$, \ie
$$F(\vec{z},\vec{w}) = \frac{1}{z_1\,\bar{w}_1\,z_2\,\bar{w}_2}.$$
Via Lemma \ref{lem:poisson_pluriharmonic}, we may determine whether this measures satisfies the Nevanlinna condition \eqref{eq:measure_Nevan}. Calculating \eg that
\begin{multline*}
    \frac{\di^2}{\di z_1 \di \bar{z}_2}\Im[z_1]\,\Im[z_2]\,F(\vec{z},\vec{z}) = -\frac{1}{4}\,\frac{\di^2}{\di z_1 \di \bar{z}_2}\,\frac{(z_1-\bar{z}_1)(z_2-\bar{z}_2)}{|z_1|^2|z_2|^2} \\
    = -\frac{1}{4} \cdot \frac{\di}{\di z_1}\,\frac{z_1-\bar{z}_1}{|z_1|^2}\cdot\frac{\di}{\di \bar{z}_2}\,\frac{z_2-\bar{z}_2}{|z_2|^2} = -\frac{1}{4}\cdot \frac{1}{{z_1}^2}\cdot\frac{-1}{{\bar{z}_2}^2} \neq 0,
\end{multline*}
sufficing to conclude that the function $(z_1,z_2) \mapsto \Im[z_1]\,\Im[z_2]\,F(\vec{z},\vec{z})$ is not pluriharmonic. \hfill$\lozenge$
\end{example}

\section{Characterization via positive semi-definite functions}
\label{sec:main_thm}

Let us begin by shortly recalling the how \HN functions in one variable are characterized via positive semi-definite functions. When $n=1$, Theorem \ref{thm:intRep_Nvar} states that $q\colon \C^+ \to \C$ is a \HN function if and only if there exists number $a \in \R$ and $b \geq 0$ as well as a positive Borel measure $\mu$ on $\R$ satisfying $\int_\R(1+t^2)^{-1}\diff\mu(t) < \infty$ such that
$$q(z) = a + b\,z + \frac{1}{\pi}\int_\R\left(\frac{1}{t-z}-\frac{t}{1+t^2}\right)\diff\mu(t)$$
for all $z \in \C^+$. This representation implies immediately that
\begin{equation}
    \label{eq:HN_decomposition_1var}
    q(z) - \bar{q(w)} = b\,(z-\bar{w}) + (z-\bar{w})\,\frac{1}{\pi}\int_\R\frac{1}{(t-z)(t-\bar{w})}\diff\mu(t)
\end{equation}
for all $z,w \in \C^+$. Equivalently, one may reformulate the above equality to say that for every \HN function $q$ the function
$$(z,w) \mapsto \frac{q(z) - \bar{q(w)}}{z-\bar{w}} = b + \frac{1}{\pi}\int_\R\frac{1}{(t-z)(t-\bar{w})}\diff\mu(t)$$
is positive semi-definite on $\C^+ \times \C^+$. Conversely, if $q\colon \C^+ \to \C$ is a holomorphic function for which the function
$$D\colon (z,w) \mapsto \frac{q(z) - \bar{q(w)}}{z-\bar{w}}$$
is positive semi-definite on $\C^+ \times \C^+$, then $q$ must be a \HN function. This is seen evaluating the function $D$ at $z=w$ and using Lemma \ref{lem:positive_semi_def_properties}. This characterization also leads to the introduction of the Nevanlinna kernel and generalized Nevanlinna functions, which we will return to in Section \ref{subsec:Nevanlinna_kernel}.

The next objective is therefore to determine whether a decomposition analogous to \eqref{eq:HN_decomposition_1var} holds for \HN functions of several variables. 

\subsection{The main theorem}
\label{subsec:main_thm}

Our main result is the following.

\begin{thm}\label{thm:positive_semi_decomposition}
Let $n \in \N$ and let $q\colon \C^{+n} \to \C$ be a holomorphic function. Then, $q$ is a \HN function if and only if there exists a vector $\vec{d} \in [0,\infty)^n$ and a positive semi-definite function $D$ on $\C^{+n} \times \C^{+n}$ satisfying the growth condition \eqref{eq:positive_semi_growth} such that the equality
\begin{equation}
    \label{eq:positive_semi_decomposition}
    q(\vec{z}) - \bar{q(\vec{w})} = \sum_{j=1}^nd_j(z_j - \bar{w}_j) + \frac{1}{(2\,\I)^{n-1}}\prod_{j=1}^n(z_j - \bar{w}_j)\cdot D(\vec{z},\vec{w})
\end{equation}
holds for all $\vec{z},\vec{w} \in \C^{+n}$.

Furthermore, let $(a,\vec{b},\mu)$ be the representing parameters of the function $q$ in the sense of Theorem \ref{thm:intRep_Nvar}. If $a = 0$, the correspondence between the function $q$ and the parameters $\vec{d}$ and $D$ is unique and it holds that $\vec{d} = \vec{b}$ and that $D$ is the Poisson-type function given by the measure $\mu$.
\end{thm}

\proof
\textit{PART 1:} Assume first that $q$ is a \HN function. Then, we wish to deduce that it admits a decomposition of the form \eqref{eq:positive_semi_decomposition} for some vector $\vec{d}$ and some positive semi-definite function $D$ on $\C^{+n} \times \C^{+n}$. Since a \HN function $q$ is uniquely determined by its data $(a,\vec{b},\mu)$ in the sense of Theorem \ref{thm:intRep_Nvar}, it may be uniquely written as the sum
$$q = q_a + q_b + q_c,$$
where $q_a$ is given by the data $(a,\vec{0},0)$, $q_b$ is given by the data $(0,\vec{b},0)$ and $q_c$ is given by the data $(0,\vec{0},\mu)$. Hence, it suffices to prove the desired result for each of these three special cases separately.

\textit{Case 1.a:} If a \HN function $q$ is given by the data $(a,\vec{0},0)$ in the sense of Theorem \ref{thm:intRep_Nvar}, then it holds, for every $\vec{z},\vec{w} \in \C^{+n}$ that $q(\vec{z}) - \bar{q(\vec{w})} = 0$. Thus, to satisfy equality \eqref{eq:positive_semi_decomposition}, one may chose $\vec{d} = \vec{0}$ and $D \equiv 0$.

\textit{Case 1.b:} If a \HN function $q$ is given by the data $(0,\vec{b},0)$ in the sense of Theorem \ref{thm:intRep_Nvar}, then it holds, for every $\vec{z},\vec{w} \in \C^{+n}$, that
$$q(\vec{z}) - \bar{q(\vec{w})} = \sum_{j=1}^nb_j(z_j - \bar{w}_j)$$
Thus, to satisfy equality \eqref{eq:positive_semi_decomposition}, one may chose $\vec{d} = \vec{b}$ and $D \equiv 0$.

\textit{Case 1.c:} If a \HN function $q$ is given by the data $(0,\vec{0},\mu)$ in the sense of Theorem \ref{thm:intRep_Nvar}, it holds, for every $\vec{z},\vec{w} \in \C^{+n}$, that
$$q(\vec{z}) - \bar{q(\vec{w})} = \frac{1}{\pi^n}\int_{\R^n}\big(K_n(\vec{z},\vec{t}) - \bar{K_n(\vec{w},\vec{t})}\big)\diff\mu(\vec{t}).$$
In order to be able to describe the difference $K_n(\vec{z},\vec{t}) - \bar{K_n(\vec{w},\vec{t})}$ more precisely, we introduce the expression
\begin{multline}\label{eq:extended_Poisson}
\pois_n(\vec{z},\vec{w},\vec{t}) := \prod_{\ell=1}^n\bigg(N_{-1}(z_\ell,t_\ell) + N_{0}(\I,t_\ell) + N_{1}(w_\ell,t_\ell)\bigg) \\ = \frac{1}{(2\,\I)^n}\prod_{\ell=1}^n\left(\frac{1}{t_\ell-z_\ell}-\frac{1}{t_\ell-\bar{w_\ell}}\right) = \frac{1}{(2\,\I)^n}\prod_{\ell=1}^n\frac{z_\ell-\bar{w_\ell}}{(t_\ell - z_\ell)(t_\ell-\bar{w_\ell})}.
\end{multline}
This expression can be thought of as an "extended" Poisson kernel of $\C^{+n}$, as $\pois_n(\vec{z},\vec{z},\vec{t})$ equals the usual Poisson kernel of $\C^{+n}$.

We now claim that the equality
\begin{equation}\label{eq:kernel_difference}
\frac{1}{2\,\I}\big(K_n(\vec{z},\vec{t}) - \bar{K_n(\vec{w},\vec{t})}\big) = \pois_n(\vec{z},\vec{w},\vec{t}) - \sum_{\substack{\vec{\rho} \in \{-1,0,1\}^n \\ -1\in\vec{\rho} \wedge 1\in\vec{\rho}}}\prod_{j=1}^nN_{\rho_j}(\epsilon_{\rho_j}(z_j,w_j),t_j)
\end{equation}
holds for every $\vec{t} \in \R^n$ and every $\vec{z},\vec{w} \in \cutN$, where the choice of variable $\epsilon_{\ell}$ is determined by
\begin{equation}\label{eq:epsilon_choice}
    \epsilon_{\ell}(\alpha,\beta) := 
\left\{\begin{array}{rcl}
\alpha & ; & \ell = -1, \\
\I & ; & \ell = 0, \\
\beta & ; & \ell = 1,
\end{array}\right.
\end{equation}
\ie the choice $\epsilon_\ell$ ensures that the the first input of the term $N_{-1}$ is always taken form the vector $\vec{z}$ and that the first input of the term $N_{1}$ is always taken form the vector $\vec{w}$. Note that this is a stronger statement than needed, as for our current goal it would suffice to consider $\vec{z},\vec{w} \in \C^{+n}$.

The proof of equality \eqref{eq:kernel_difference} follows closely the proof of the special case when $\vec{z} = \vec{w}$, presented in \cite[Prop. 3.3]{LugerNedic2019}. To that end, we observe that
$$\bar{K_n(\vec{w},\vec{t})} = -\I\bigg(2\prod_{j=1}^n\big(N_{1}(w_j,t_j) +  N_{0}(\I,t_j)\big)-\prod_{j=1}^nN_{0}(\I,t_j)\bigg).$$
Expanding now the following products as sums yields
$$\begin{array}{RCL}
\multicolumn{3}{L}{\prod_{j=1}^n\big(N_{-1}(z_j,t_j)+N_{0}(\I,t_j)+N_{1}(w_j,t_j)\big)} \\[0.45cm]
~ & = & \sum_{\vec{\rho} \in \{-1,0,1\}^n}\prod_{j=1}^nN_{\rho_j}(\epsilon_{\rho_j}(z_j,w_j),t_j), \\[0.45cm]
\prod_{j=1}^n\big(N_{1}(w_j,t_j)+N_{0}(\I,t_j)\big) & = & \sum_{\substack{\vec{\rho} \in \{-1,0,1\}^n \\ -1\not\in\vec{\rho}}}\prod_{j=1}^nN_{\rho_j}(\epsilon_{\rho_j}(z_j,w_j),t_j), \\[0.45cm]
\prod_{j=1}^n\big(N_{-1}(z_j,t_j)+N_{0}(\I,t_j)\big) & = & \sum_{\substack{\vec{\rho} \in \{-1,0,1\}^n \\ 1\not\in\vec{\rho}}}\prod_{j=1}^nN_{\rho_j}(\epsilon_{\rho_j}(z_j,w_j),t_j), \\[0.45cm]
\prod_{j=1}^nN_{0}(\I,t_j) & = & \sum_{\substack{\vec{\rho} \in \{-1,0,1\}^n \\ -1\not\in\vec{\rho} \wedge 1\not\in\vec{\rho}}}\prod_{j=1}^nN_{\rho_j}(\epsilon_{\rho_j}(z_j,w_j),t_j),
\end{array}$$
and hence 
$$\begin{array}{RCL}
\multicolumn{3}{L}{\frac{1}{2\,\I}\big(K_n(\vec{z},\vec{t}) - \bar{K_n(\vec{w},\vec{t})}\big)} \\[0.45cm]
~ & = & \prod_{j=1}^n\big(N_{-1}(z_j,t_j)+N_{0}(\I,t_j)\big) + \prod_{j=1}^n\big(N_{1}(w_j,t_j)+N_{0}(\I,t_j)\big) - \prod_{j=1}^nN_{0}(\I,t_j) \\[0.45cm]
~ & = & \sum_{\vec{\rho} \in \{-1,0,1\}^n}\prod_{j=1}^nN_{\rho_j}(\epsilon_{\rho_j}(z_j,w_j),t_j)-\sum_{\substack{\vec{\rho} \in \{-1,0,1\}^n \\ -1\in\vec{\rho} \wedge 1\in\vec{\rho}}}\prod_{j=1}^nN_{\rho_j}(\epsilon_{\rho_j}(z_j,w_j),t_j) \\[0.45cm]
~ & = & \pois_n(\vec{z},\vec{w},\vec{t}) - \sum_{\substack{\vec{\rho} \in \{-1,0,1\}^n \\ -1\in\vec{\rho} \wedge 1\in\vec{\rho}}}\prod_{j=1}^nN_{\rho_j}(\epsilon_{\rho_j}(z_j,w_j),t_j),
\end{array}$$
as desired.

Using formula \eqref{eq:kernel_difference}, we deduce that
\begin{multline}\label{eq:proof_symmetric_temp1}
q(\vec{z}) - \bar{q(\vec{w})} \\ = \frac{2\,\I}{\pi^n}\int_{\R^n} \pois_n(\vec{z},\vec{w},\vec{t})\diff\mu(\vec{t}) - \sum_{\substack{\vec{\rho} \in \{-1,0,1\}^n \\ -1\in\vec{\rho} \wedge 1\in\vec{\rho}}}\frac{2\,\I}{\pi^n}\int_{\R^n}\prod_{j=1}^nN_{\rho_j}(\epsilon_{\rho_j}(z_j,w_j),t_j)\diff\mu(\vec{t}).
\end{multline}
Since
\begin{multline*}
\frac{2\,\I}{\pi^n}\int_{\R^n} \pois_n(\vec{z},\vec{w},\vec{t})\diff\mu(\vec{t}) \\ = \frac{1}{(2\,\I)^{n-1}}\prod_{j=1}^n(z_j - \bar{w}_j) \cdot \frac{1}{\pi^n}\int_{\R^n}\prod_{j=1}^n\frac{1}{(t_j-z_j)(t_j-\bar{w}_j)}\diff\mu(\vec{t}),    
\end{multline*}
due to the definition of $\pois_n$, we may choose $D$ as the Poisson-type function given by the measure $\mu$. This is indeed a valid choice, as a Poisson-type function is positive semi-definite on $\C^{+n} \times \C^{+n}$ due to Lemma \ref{lem:positive_semi_poisson} and satisfies the growth condition \eqref{eq:positive_semi_growth} due to Lemma \ref{lem:poisson_growth}.  

Finally, we claim that
$$\int_{\R^n}\prod_{j=1}^nN_{\rho_j}(\epsilon_{\rho_j}(z_j,w_j),t_j)\diff\mu(\vec{t}) = 0$$
for every vectors $\vec{z},\vec{w} \in \C^{+n}$ and every indexing vector $\vec{\rho}$ as in the sum in formula \eqref{eq:proof_symmetric_temp1}, \ie with at lest one entry equal to $1$ and at least one entry equal to $-1$. To infer this, we observe that once an indexing vector $\vec{\rho}$ has been chosen, we may combine the vectors $\vec{z},\vec{w} \in \C^{+n}$ into a single vector $\vec{\xi} \in \C^{+n}$ via setting
$$\xi_j := \left\{\begin{array}{rcl}
z_j & ; & \rho_j = -1, \\
\I & ; & \rho_j = 0, \\
w_j & ; & \rho_j = 1.
\end{array}\right.$$
Hence,
$$N_{\rho_j}(\epsilon_{\rho_j}(z_j,w_j),t_j) = N_{\rho_j}(\xi_j,t_j)$$
and the desired result follows from the fact that the measure $\mu$ satisfies the Nevnalinna condition \eqref{eq:measure_Nevan}, \cf \cite[Thm. 5.1]{LugerNedic2019}. Note that it is important that we may write the choice $\epsilon_{\rho_j}$ in terms of a single vector from $\C^{+n}$, as \cite[Thm. 5.1 (b)]{LugerNedic2019} only implies the desired result in the case where the first inputs of the terms $N_{-1}$ and $N_{1}$ are taken from the same vector. 

In conclusion, formula \eqref{eq:proof_symmetric_temp1} provides a decomposition of the form \eqref{eq:positive_semi_decomposition} where the function $D$ is chosen as described above and $\vec{d} = \vec{0}$. This finishes part 1 of the proof.

\textit{PART 2:} Assume now that $q\colon \C^{+n} \to \C$ is a holomorphic function for which there exists some vector $\vec{d} \in [0,\infty)^n$ and some positive semi-definite function $D$ on $\C^{+n} \times \C^{+n}$ such that equality \eqref{eq:positive_semi_decomposition} holds for all $\vec{z},\vec{w} \in \C^{+n}$. To show that $q$ must be a \HN function, we only need to check that $\Im[q] \geq 0$. To that end, we may choose $\vec{z} = \vec{w}$ in equality \eqref{eq:positive_semi_decomposition} to get
$$2\,\I\,\Im[q(\vec{z})] = 2\,\I\sum_{j=1}^nd_j\Im[z_j] + 2\,\I\,\prod_{j=1}^n\Im[z_j]\cdot D(\vec{z},\vec{z}).$$
Diving both sides by $2\,\I$ and noting that $D(\vec{z},\vec{z}) \geq 0$ by Lemma \ref{lem:positive_semi_def_properties} yields the desired result.

As such, the only remaining step is to show that the vector $\vec{d}$ must be equal to $\vec{b}$ and that the function $D$ must be the Poisson-type function given by the measure $\mu$. For the former, we calculate using formula \eqref{eq:b_parameters} that, one one hand,
$$\lim\limits_{z_\ell \ntto \infty}\frac{q(\vec{z}) - \bar{q(\vec{w})}}{z_\ell} = b_\ell.$$
On the other hand, using decomposition \eqref{eq:positive_semi_decomposition}, we have
$$\lim\limits_{z_\ell \ntto \infty}\frac{q(\vec{z}) - \bar{q(\vec{w})}}{z_\ell} = d_\ell + \frac{1}{(2\,\I)^{n-1}}\prod_{\substack{j = 1 \\ j \neq \ell}}^n(z_j - \bar{w}_j)\cdot\lim\limits_{z_\ell \ntto \infty}D(\vec{z},\vec{w}) = d_\ell$$
due to the assumption that the function $D$ satisfies the growth condition \eqref{eq:positive_semi_growth}.

Since we know by now that the function $q$ is a \HN function, it is represented by some data $(a,\vec{b},\mu)$ in the sense of Theorem \ref{thm:intRep_Nvar}. By Part 1 of the proof, it holds that
$$q(\vec{z}) - \bar{q(\vec{w})} = \sum_{j=1}^nb_j(z_j - \bar{w}_j) + q_c(\vec{z}) - \bar{q_c(\vec{w})},$$
where the function $q_c$ is a \HN function given by the data $(0,0,\mu)$. However, by assumption, the function $q$ also admits a decomposition of the form \eqref{eq:positive_semi_decomposition}, where we already know that $b_\ell = d_\ell$ for all $\ell \in \{1,\ldots,n\}$. Comparing this decomposition with the one above yields that
$$D(\vec{z},\vec{w}) = \frac{(2\,\I)^{n-1}}{\prod_{j=1}^n(z_j - \bar{w}_j)}\big(q_c(\vec{z}) - \bar{q_c(\vec{w})}\big).$$
However, by Step 1.c of Part 1 of the proof, it holds that
$$q_c(\vec{z}) - \bar{q_c(\vec{w})} = \frac{1}{(2\,\I)^{n-1}}\prod_{j=1}^n(z_j - \bar{w}_j)\cdot\frac{1}{\pi^n}\int_{\R^n}\prod_{\ell=1}^{n}\frac{1}{(t_\ell - z_\ell)(t_\ell - \bar{w}_\ell)}\diff\mu(\vec{t}),$$
implying that the function $D$ is necessarily the Poisson-type function given by the measure $\mu$, finishing Part 2 of the proof.
\endproof

A slight reformulation of this result can be stated as follows.

\begin{coro}
\label{coro:D_rep}
A function $q\colon \C^{+n} \to \C$ is a \HN function if and only if there exists a number $a \in \R$, a vector $\vec{b} \in [0,\infty)^n$ and a Poisson-type function $D$ satisfying condition \eqref{eq:poisson_pluriharmonic} such that the formula
\begin{equation}
    \label{eq:D_rep}
    q(\vec{z}) = \big(a - \I\,D(\I\,\vec{1},\I\,\vec{1})\big) + \sum_{j=1}^nb_j\,z_j + \frac{1}{(2\,\I)^{n-1}}\prod_{j=1}^n(z_j + \I) \cdot D(\vec{z},\I\,\vec{1})
\end{equation}
holds for all $\vec{z} \in \C^{+n}$.
\end{coro}

\proof
Note first that the integral of the kernel $K_n$ with respect to a measure $\mu$ satisfying the growth condition may be written in terms of the Poisson-type function $D$ given by the same measure as
$$\frac{1}{\pi^n}\int_{\R^n}K_n(\vec{z},\vec{t})\diff\mu(\vec{t}) = \frac{1}{(2\,\I)^{n-1}}\prod_{j=1}^n(z_j + \I) \cdot D(\vec{z},\I\,\vec{1}) - \I\,D(\I\,\vec{1},\I\,\vec{1}).$$
The statement of the corollary now follows from Theorem \ref{thm:intRep_Nvar}, Theorem \ref{thm:positive_semi_decomposition} and Lemma \ref{lem:poisson_pluriharmonic}.
\endproof

\subsection{The Nevanlinna kernel in several variables}
\label{subsec:Nevanlinna_kernel}

For any holomorphic function $q\colon \C^+ \to \C$, one can consider its \emph{Nevanlinna kernel} $\kilb_q\colon \C^+ \times \C^+ \to \C$, defined by
$$\mathcal{K}_{q}(z,w) := \frac{q(z) - \bar{q(w)}}{z-\bar{w}}.$$
In general, this kernel may be considered with regards to either a scalar-, matrix- or operator-valued function $q$, see \eg \cite{DahoLanger1985,KreinLanger1973,KreinLanger1977,KreinLanger1981}. As summarized in the introduction to Section \ref{sec:main_thm}, it holds that a function $q\colon \C^+ \to \C$ is a \HN function if and only if $\kilb_q$ is positive semi-definite.

An analogous characterization for \HN functions of several variables, based on Theorem \ref{thm:positive_semi_decomposition}, is the following.

\begin{thm}
\label{thm:Nevan_kernel_Nvar}
Let $q$ be a holomorphic function on $\C^{+n}$. Then, $q$ is a \HN function if and only if the function
$$(\vec{z},\vec{w}) \mapsto \frac{(2\,\I)^{n-1}}{\prod_{j=1}^n(z_j - \bar{w}_j)}\big(q(\vec{z}) - \bar{q(\vec{w})}\big)$$
is positive semi-definite on $\C^{+n} \times \C^{+n}$. In that case, it holds that
$$\frac{(2\,\I)^{n-1}}{\prod_{j=1}^n(z_j - \bar{w}_j)}\big(q(\vec{z}) - \bar{q(\vec{w})}\big) = \sum_{j=1}^nb_j\prod_{\substack{\ell = 1 \\ \ell \neq j}}^n\frac{2\,\I}{z_\ell - \bar{w}_\ell} + D(\vec{z},\vec{w}),$$
where the vector $\vec{b}$ is as in representation \eqref{eq:intRep_Nvar} and the function $D$ is as described in Theorem \ref{thm:positive_semi_decomposition}.
\end{thm}

\proof
If the function $q$ is represented by the data $(a,\vec{0},0)$ or $(0,\vec{0},\mu)$, then the result follows immediately by steps 1.a and 1.c of the proof of Theorem \ref{thm:positive_semi_decomposition}. If, instead, the function $q$ is represented by the data $(0,\vec{b},0)$, it holds that
$$\frac{(2\,\I)^{n-1}}{\prod_{j=1}^n(z_j - \bar{w}_j)}\big(q(\vec{z}) - \bar{q(\vec{w})}\big) = \sum_{j=1}^nb_j\prod_{\substack{\ell = 1 \\ \ell \neq j}}^n\frac{2\,\I}{z_\ell - \bar{w}_\ell}.$$
Since the function
$$(z,w) \mapsto \frac{2\,\I}{z-\bar{w}}$$
is positive semi-definite on $\C^+ \times \C^+$, the result then follows by Lemma \ref{lem:positive_semi_def_properties}.

For the converse statement, we only need to show that $\Im[q(\vec{z})] \geq 0$ which may be done as in Part 2 of the proof of Theorem \ref{thm:positive_semi_decomposition}.
\endproof

\begin{coro}
Let $F$ be a positive semi-definite function on $\C^{+n} \times \C^{+n}$. Then,
\begin{equation}
\label{eq:F_till_q}
    F(\vec{z},\vec{w}) = \frac{(2\,\I)^{n-1}}{\prod_{j=1}^n(z_j - \bar{w}_j)}\big(q(\vec{z}) - \bar{q(\vec{w})}\big)
\end{equation}
for some \HN function $q$ if and only if there exists a vector $\vec{d} \in [0,\infty)^n$ and a Poisson-type function $D$ satisfying condition \eqref{eq:poisson_pluriharmonic} such that
\begin{equation}
    \label{eq:D_till_F}
    F(\vec{z},\vec{w}) = \sum_{j=1}^nd_j\prod_{\substack{\ell = 1 \\ \ell \neq j}}^n\frac{2\,\I}{z_\ell - \bar{w}_\ell} + D(\vec{z},\vec{w})
\end{equation}
for all $\vec{z},\vec{w} \in \C^{+n}$.
\end{coro}

\proof
If the function $F$ can be written in terms of a \HN function $q$ as in formula \eqref{eq:F_till_q}, then the result follows by Theorem \ref{thm:positive_semi_decomposition} and Theorem \ref{thm:Nevan_kernel_Nvar}.

Conversely, assume that there exists a vector $\vec{d}$ and a Poisson-type function $D$ as in the Corollary such that equality \eqref{eq:D_till_F} holds. Let $\sigma$ be the positive Borel measure satisfying the growth condition \eqref{eq:measure_growth} that gives the function $D$. By Lemma \ref{lem:poisson_pluriharmonic}, the measure sigma also satisfies the Nevanlinna condition \eqref{eq:measure_Nevan} due to the assumption on the function $D$. Hence, we may define a \HN function $q$ via Theorem \ref{thm:intRep_Nvar} using the data $(0,\vec{d},\sigma)$. The result now follows by Theorem \ref{thm:Nevan_kernel_Nvar}.
\endproof

Using Theorem \ref{thm:Nevan_kernel_Nvar}, we may now propose a multidimensional analogue to the classical Nevaninna kernel that generalizes its core property of characterizing \HN functions.

\begin{define}
\label{def:Nevan_kernel_Nvar}
Let $q\colon\C^{+n} \to \C$ be a holomorphic function. Then, its \emph{Nevanlinna kernel} $\kilb_q\colon \C^{+n} \times \C^{+n} \to \C$ is defined as
$$\mathcal{K}_{q}(\vec{z},\vec{w}) := \frac{(2\,\I)^{n-1}}{\prod_{j=1}^n(z_j - \bar{w}_j)}\big(q(\vec{z}) - \bar{q(\vec{w})}\big).$$
\end{define}

When $n=1$, the Nevanlinna kernel is used to define \emph{generalized Nevanlinna functions}. A meromorphic function $q\colon \C^+ \to \C$ with domain of holomorphy $\Dom(q) \subseteq \C^+$ is called a generalized Nevanlinna function of class $\mathcal{N}_\kappa(\C^+)$ if its Nevanlinna kernel $\kilb_q$ has $\kappa$ negative squares \cite[pg. 187]{KreinLanger1977}. We recall that $\mathcal{K}_{q}$ having $\kappa$ negative squares means that for arbitrary $N \in \N$ and arbitrary $z_1,\ldots,z_N \in \C^+$ the matrix
$$\left(\mathcal{K}_{q}(z_i,z_j)\right)_{i,j = 1}^N$$
has at most $\kappa$ negative eigenvalues and $\kappa$ is minimal with this property.

Using Definition \ref{def:Nevan_kernel_Nvar}, generalized Nevanlinna functions of several variables may be introduced completely analogously. 

\begin{define}
A meromorphic function $q\colon \C^{+n} \to \C$ with domain of holomorphy $\Dom(q) \subseteq \C^{+n}$ is called a generalized Nevanlinna function of class $\mathcal{N}_\kappa(\C^{+n})$ if its Nevanlinna kernel $\kilb_q$ has $\kappa$ negative squares.
\end{define}

The detailed study of this class of functions lies outside the scope of this paper.

\subsection{Decomposition of the symmetric extension}
\label{subsec:symmetric_extension}

We recall that the integral representation in formula \eqref{eq:intRep_Nvar} is well-defined for any $\vec{z} \in \cutN$, which may be used to extend any Herglotz-Nevanlinna function $q$ from $\C^{+n}$ to $\cutN$. This extension is called the \emph{symmetric extension} of the function $q$ and is denoted as $q_\mathrm{sym}$. We note that that the symmetric extension of a \HN function $q$ is different from its possible analytic extension as soon as $\mu \neq 0$, \cf \cite[Prop. 6.10]{LugerNedic2019}. The symmetric extension also satisfies the following variable-dependence property, \cf \cite[Prop. 6.9]{LugerNedic2019}.

Just as a \HN function $q$ can always be symmetrically extended to $\cutN$ via its integral representation \eqref{eq:intRep_Nvar}, so too can we consider the symmetric extension of a Poisson-type function. This symmetric extension will automatically be positive semi-definite on $\cutN \times \cutN$ as the proof of Lemma \ref{lem:positive_semi_poisson} still remains valid even if the variables $\vec{z}$ and $\vec{w}$ are taken from $\cutN$ instead. However, a direct analogue of Theorem \ref{thm:positive_semi_decomposition} does not hold, as the following example shows.

\begin{example}\label{ex:standard}
Consider the \HN function $q$ given by $q(z_1,z_2) = -(z_1+z_2)^{-1}$ for $(z_1,z_2) \in \C^{+2}$. This function in represented by the data $(0,\vec{0},\mu)$ in the sense of Theorem \ref{thm:intRep_Nvar}, where the measure $\mu$ is defined, for any Borel set $U \subseteq \R^2$, as
$$\mu(U) := \pi\int_\R\chi_U(t,-t)\diff t.$$
Here, $\chi_U$ denotes the characteristic function of the set $U$. The function $D$ from Theorem \ref{thm:positive_semi_decomposition} can then be calculated using standard residue calculus and equals
$$D(\vec{z},\vec{w}) = \frac{2\,\I\,(z_1+z_2-\bar{w}_1 -\bar{w}_2)}{(z_1-\bar{w}_1)(z_2-\bar{w}_2)(z_1+z_2)(\bar{w}_1+\bar{w}_2)}.$$

Furthermore, also using standard residue calculus, the symmetric extension of the function $q$ can be calculated to be
$$q_\mathrm{sym}(z_1,z_2) = \left\{\begin{array}{rcl}
-(z_1+z_2)^{-1} & ; & (z_1,z_2) \in \C^+ \times \C^+, \\
(\I - z_1)^{-1} & ; & (z_1,z_2) \in \C^- \times \C^+, \\
(\I - z_2)^{-1} & ; & (z_1,z_2) \in \C^+ \times \C^-, \\
(z_1+z_2)^{-1}+(\I - z_1)^{-1}+(\I - z_2)^{-1} & ; & (z_1,z_2) \in \C^- \times \C^-,
\end{array}\right.$$
while the values of $D_\mathrm{sym}$ are shown in Table \ref{tab:symetric_extension}.

\begin{table}[!ht]
    \centering
    \small{$$\begin{array}{C|C|C}
        D_\mathrm{sym}(\vec{z},\vec{w}) & \vec{z} \in & \vec{w} \in  \\
        \hline
        \frac{2\,\I\,(z_1+z_2-\bar{w}_1 -\bar{w}_2)}{(z_1-\bar{w}_1)(z_2-\bar{w}_2)(z_1+z_2)(\bar{w}_1+\bar{w}_2)} & \C^+ \times \C^+ & \C^+ \times \C^+ \\[0.4cm]
        \frac{2\,\I}{(z_1+\bar{w}_2)(z_2-\bar{w}_2)(\bar{w}_1+\bar{w}_2)} & \C^- \times \C^+ & \C^+ \times \C^+ \\[0.4cm]
        \frac{2\,\I}{(z_2+\bar{w}_1)(z_1-\bar{w}_1)(\bar{w}_1+\bar{w}_2)} & \C^+ \times \C^- & \C^+ \times \C^+ \\[0.4cm]
        \frac{2\,\I\,(z_1+z_2+\bar{w}_1 +\bar{w}_2)}{(z_1+\bar{w}_2)(z_2+\bar{w}_1)(z_1+z_2)(\bar{w}_1+\bar{w}_2)} & \C^- \times \C^- & \C^+ \times \C^+ \\[0.4cm]
        \hline
        \frac{2\,\I}{(z_1 + z_2)(z_2 + \bar{w}_1)(z_2 - \bar{w}_2)} & \C^+ \times \C^+ & \C^- \times \C^+ \\[0.4cm]
        \frac{2\,\I\,(z_1-z_2-\bar{w}_1 +\bar{w}_2)}{(z_1-\bar{w}_1)(z_2+\bar{w}_1)(z_1 + \bar{w}_2)(z_2 - \bar{w}_2)} & \C^- \times \C^+ & \C^- \times \C^+ \\[0.4cm]
        0 & \C^+ \times \C^- & \C^- \times \C^+ \\[0.4cm]
        \frac{-2\,\I}{(z_1 + z_2)(z_1 - \bar{w}_1)(z_1 + \bar{w}_2)} & \C^- \times \C^- & \C^- \times \C^+ \\[0.4cm]
        \hline
        \frac{2\,\I}{(z_1 + z_2)(z_1 - \bar{w}_1)(z_1 + \bar{w}_2)} & \C^+ \times \C^+ & \C^+ \times \C^- \\[0.4cm]
        0 & \C^- \times \C^+ & \C^+ \times \C^- \\[0.4cm]
        \frac{2\,\I\,(-z_1+z_2+\bar{w}_1-\bar{w}_2)}{(z_1-\bar{w}_1)(z_2+\bar{w}_1)(z_1 + \bar{w}_2)(z_2 - \bar{w}_2)} & \C^+ \times \C^- & \C^+ \times \C^- \\[0.4cm]
        \frac{-2\,\I}{(z_1 + z_2)(z_2 + \bar{w}_1)(z_2 - \bar{w}_2)} & \C^- \times \C^- & \C^+ \times \C^- \\[0.4cm]
        \hline
        \frac{-2\,\I\,(z_1+z_2+\bar{w}_1 +\bar{w}_2)}{(z_1+\bar{w}_2)(z_2+\bar{w}_1)(z_1+z_2)(\bar{w}_1+\bar{w}_2)} & \C^+ \times \C^+ & \C^- \times \C^- \\[0.4cm]
        \frac{-2\I}{(z_1 - \bar{w}_1)(z_2 + \bar{w}_1)(\bar{w}_1 + \bar{w}_2)} & \C^- \times \C^+ & \C^- \times \C^- \\[0.4cm]
        \frac{-2\I}{(z_2 - \bar{w}_2)(z_1 + \bar{w}_2)(\bar{w}_1 + \bar{w}_2)} & \C^+ \times \C^- & \C^- \times \C^- \\[0.4cm]
        \frac{2\I(-z_1 - z_2 + \bar{w}_1 + \bar{w}_2)}{(z_1 + z_2)(z_1 - \bar{w}_1)(z_2 - \bar{w}_2)(\bar{w}_1 + \bar{w}_2)} & \C^- \times \C^- & \C^- \times \C^-
    \end{array}$$}
    \caption{The symmetric extension of the function $D$ from Example \ref{ex:standard}}
    \label{tab:symetric_extension}
\end{table}

If we now choose \eg $\vec{z} \in \C^+ \times \C^-$ and $\vec{w} \in \C^- \times \C^+$, we see that
$$\frac{1}{\I - \bar{z}_2} + \frac{1}{\I+\bar{w}_1} = q_\mathrm{sym}(\vec{z}) - \bar{q_\mathrm{sym}(\vec{w})} \neq \frac{1}{(2\,\I)^{n-1}}\prod_{j=1}^n(z_j - \bar{w}_j)\cdot D_\mathrm{sym}(\vec{z},\vec{w}) = 0,$$
implying the presence of an error term. \hfill$\lozenge$
\end{example}

The following proposition gives a characterization of the symmetric extension of a \HN functions in analogy to Theorem \ref{thm:positive_semi_decomposition}.

\begin{prop}\label{prop:symmetric_positive_semi_decomposition}
Let $n \in \N$ and let $f\colon \cutN \to \C$ be a holomorphic function. Then, $f = q_\mathrm{sym}$ for some \HN function $q$ if and only if there exists a vector $\vec{d} \in [0,\infty)^n$ and a positive Borel measure $\sigma$ on $\R^n$ satisfying the growth condition \eqref{eq:measure_growth} and the Nevanlinna condition \eqref{eq:measure_Nevan} such that the equality
\begin{multline}
    \label{eq:symmetric_positive_semi_decomposition}
    f(\vec{z}) - \bar{f(\vec{w})} \\ = \sum_{j=1}^nd_j(z_j - \bar{w}_j) + \frac{1}{(2\,\I)^{n-1}}\prod_{j=1}^n(z_j - \bar{w}_j)\cdot D_\mathrm{sym}(\vec{z},\vec{w}) - E(\vec{z},\vec{w})
\end{multline}
holds for all $\vec{z},\vec{w} \in \cutN$. Here, $D_\mathrm{sym}$ denotes the symmetric extension of the Poisson-type function $D$ given by the measure $\sigma$, the error term $E$ is defined as
$$E(\vec{z},\vec{w}) := \sum_{\substack{\vec{\rho} \in \{-1,0,1\}^n \\ -1\in\vec{\rho} \wedge 1\in\vec{\rho}}}\frac{2\,\I}{\pi^n}\int_{\R^n}\prod_{j=1}^nN_{\rho_j}(\epsilon_{\rho_j}(z_j,w_j),t_j)\diff\sigma(\vec{t})$$
and the choice of variable $\epsilon_\ell$ is defined by formula \eqref{eq:epsilon_choice}.
\end{prop}

\proof
Assume first that $f = q_\mathrm{sym}$ for some \HN function $q$ given by the data $(a,\vec{b},\mu)$ in the sense of Theorem \ref{thm:intRep_Nvar}. As in the proof of Theorem \ref{thm:positive_semi_decomposition}, we separately investigate three case with respect to the representing parameters of the function $q$. Cases 1.a and 1.b, \ie when the function $q$ is either represented by data of the form $(a,\vec{0},0)$ or $(0,\vec{b},0)$ may be considered completely analogously as before.

If, instead, a \HN function $q$ is given by the data $(0,\vec{0},\mu)$ in the sense of Theorem \ref{thm:intRep_Nvar}, it holds that
$$q_\mathrm{sym}(\vec{z}) - \bar{q_\mathrm{sym}(\vec{w})} = \frac{1}{\pi^n}\int_{\R^n}\big(K_n(\vec{z},\vec{t}) - \bar{K_n(\vec{w},\vec{t})}\big)\diff\mu(\vec{t})$$
for every $\vec{z},\vec{w} \in \cutN$. Using formula \eqref{eq:kernel_difference}, which we already originally proved holds for $\vec{z},\vec{w} \in \cutN$, we deduce that
\begin{multline*}
q_\mathrm{sym}(\vec{z}) - \bar{q_\mathrm{sym}(\vec{w})} \\ = \frac{2\,\I}{\pi^n}\int_{\R^n} \pois_n(\vec{z},\vec{w},\vec{t})\diff\mu(\vec{t}) - \sum_{\substack{\vec{\rho} \in \{-1,0,1\}^n \\ -1\in\vec{\rho} \wedge 1\in\vec{\rho}}}\frac{2\,\I}{\pi^n}\int_{\R^n}\prod_{j=1}^nN_{\rho_j}(\epsilon_{\rho_j}(z_j,w_j),t_j)\diff\mu(\vec{t}).
\end{multline*}
Since
$$\frac{2\,\I}{\pi^n}\int_{\R^n} \pois_n(\vec{z},\vec{w},\vec{t})\diff\mu(\vec{t}) = \frac{1}{(2\,\I)^{n-1}}\prod_{j=1}^n(z_j - \bar{w}_j)\cdot D_\mathrm{sym}(\vec{z},\vec{w}),$$
the result follows.

Conversely, define a \HN function $q$ via representation \eqref{eq:intRep_Nvar} using the data $(\Re[f(\I\,\ldots,\I)],\vec{d},\sigma)$. By what we just proved, it holds that
\begin{equation}
    \label{eq:temp1}
    f(\vec{z}) - \bar{f(\vec{w})} = q_\mathrm{sym}(\vec{z}) - \bar{q_\mathrm{sym}(\vec{w})}
\end{equation}
for all $\vec{z},\vec{w} \in \cutN$. When $\vec{z} = \vec{w}$, the above equality implies that
$$\Im[f(\vec{z})] = \Im[f(\vec{q})]$$
for all $\vec{z} \in \cutN$. Therefore,
$$f(\vec{z}) = q(\vec{z}) + C(\vec{z}),$$
where the function $C$ equals a real constant on each connected component of $\cutN$. By construction,
$$\Re[q(\I,\ldots,\I)] = \Re[f(\I,\ldots,\I)],$$
implying that $C \equiv 0$ on $\C^{+n}$ and hence $f = q_\mathrm{sym}$ on $\C^{+n}$. Finally, when $\vec{z} \in \cutN$ is arbitrary and $\vec{w} = (\I,\ldots,\I)$, equality \eqref{eq:temp1} implies that $f = q_\mathrm{sym}$ on $\cutN$. This finishes the proof.
\endproof

\subsection{Loewner functions}
\label{subsec:loewner}

Theorem \ref{thm:positive_semi_decomposition} provides a universal description of the difference $q(\vec{z}) - \bar{q(\vec{w})}$ for a \HN function $q$ and $\vec{z},\vec{w} \in \C^{+n}$. Another approach one can take is to ask instead that this difference may be written in a specific form. This leads us to consider the following class of functions, \cf \cite[pg. 3003]{AglerEtal2016}.

\begin{define}\label{def:lowner_class}
A holomorphic function $h\colon \C^{+n} \to \C$ is called a \emph{Loewner function} (on $\C^{+n}$) if there exist $n$ positive semi-definite functions $F_1,\ldots,F_n$ on $\C^{+n} \times \C^{+n}$, such that the equality
\begin{equation}
    \label{eq:lowner_functions}
    h(\vec{z}) - \bar{h(\vec{w})} = \sum_{\ell = 1}^n(z_\ell - \bar{w}_\ell)F_\ell(\vec{z},\vec{w}).
\end{equation}
holds for all $\vec{z},\vec{w} \in \C^{+n}$.
\end{define}

Due to Lemma \ref{lem:positive_semi_def_properties}, it is easily seen that every Loewner function is also a \HN function. The converse result is dependant on the number of variables $n$. When $n=1$, it follows from the results summarized in the beginning of Section \ref{sec:main_thm} that every \HN function is also a Loewner function. If $n=2$, this still holds, with the result being a consequence of a theorem concerning Schur function on the polydisk, see \cite{Agler1990} and \cite[Thm. 1.4]{AglerEtal2016}. When $n \geq 3$, Loewner functions form a proper subclass of \HN functions, with the result being a consequence of the theory of commuting contractions \cite{Parrott1970,Varopoulos1974}.

Using the Nevanlinna kernel, we may describe Loewner functions in the following way.

\begin{prop}
Let $h\colon \C^{+n} \to \C$ be a holomorphic function. Then, $h$ is a Loewner function if and only if there exist $n$ positive semi-definite functions $F_1,\ldots,F_n$ on $\C^{+n} \times \C^{+n}$ such that the equality
\begin{equation}
    \label{eq:Loewner_kernel}
    \kilb_h(\vec{z},\vec{w}) = \sum_{\ell=1}^nF_\ell(\vec{z},\vec{w})\prod_{\substack{j=1\\j\neq\ell}}^n\frac{2\,\I}{z_j-\bar{w}_j}
\end{equation}
holds for all $\vec{z},\vec{w} \in \C^{+n}$.
\end{prop}

\proof
If $h$ is a Loewner function, it admits a decomposition of the form \eqref{eq:lowner_functions} for some positive semi-definite functions $F_1,\ldots,F_n$ on $\C^{+n} \times \C^{+n}$. Using this to rewrite the difference $h(\vec{z}) - \bar{h(\vec{w})}$ in the definition of $\kilb_h$ gives equality \eqref{eq:Loewner_kernel}.

Conversely, assuming equality \eqref{eq:Loewner_kernel} holds for some positive semi-definite functions $F_1,\ldots,F_n$ on $\C^{+n} \times \C^{+n}$. Then, multiplying both sides of the equality with $(2\,\I)^{-n+1}\prod_{j=1}^n(z_j-\bar{w}_j)$ gives a decomposition of the from \eqref{eq:lowner_functions}, as desired.
\endproof

For $n=1$, decomposition \eqref{eq:lowner_functions} coincides with decomposition \eqref{eq:positive_semi_decomposition}, implying that a given Loewner function admits precisely one decomposition of the form \eqref{eq:lowner_functions}. For $n \geq 2$, decomposition \eqref{eq:lowner_functions} is not-necessarily unique and the functions $F_1,\ldots,F_n$ are not necessarily Poisson-type functions. To illustrate this, consider the following two examples.

\begin{example}
Let $h(z_1,z_2) := \I$ for all $(z_1,z_2) \in \C^{+2}$. Then, it holds that
\begin{eqnarray*}
h(z_1,z_2) - \bar{h(w_1,w_2)} & = & (z_1 - \bar{w}_1)\cdot\frac{2\,\I}{z_1-\bar{w}_1} + (z_2-\bar{w}_2)\cdot0 \\
~ & = & (z_1 - \bar{w}_1)\cdot0 + (z_2-\bar{w}_2)\cdot\frac{2\,\I}{z_2-\bar{w}_2} \\
~ & = & (z_1 - \bar{w}_1)\cdot\frac{2k_1\I}{z_1-\bar{w}_1} + (z_2-\bar{w}_2)\cdot\frac{2k_2\I}{z_2-\bar{w}_2},
\end{eqnarray*}
where $k_1,k_2 \geq 0$ with $k_1+k_2 = 1$. Thus, this function $h$ admits infinitely many different decompositions of the form \eqref{eq:lowner_functions}. It is also noteworthy that while the function
$$(z_1,w_1) \mapsto \frac{2\,\I}{z_1-\bar{w}_1}$$
is a Poisson-type function on $\C^+ \times \C^+$, the function
$$F\colon ((z_1,z_2),(w_1,w_2)) \mapsto \frac{2\,\I}{z_1-\bar{w}_1}$$
is not a Poisson-type function on $\C^{+2} \times \C^{+2}$, though it still is positive semi-definite on $\C^{+2} \times \C^{+2}$. Indeed, if it were, it would be given by some measure $\mu$, which one would be able to reconstruct via the Stieltjes inversion formula \eqref{eq:Stieltjes_for_Poisson}. Taking $\psi$ to be any function as in Lemma \ref{lem:Stieltjes_for_Poisson}, we calculate that
\begin{multline*}
    0 \leq \iint_{\R^2}|\psi(x_1,x_2)|\,y_1\,y_2\,F(\vec{x}+\I\,\vec{y},\vec{x}+\I\,\vec{y})\diff\vec{x} = \iint_{\R^2}|\psi(x_1,x_2)|\,y_2\diff\vec{x} \\
    \leq  C\,\pi^2\,y_2 \xrightarrow{\vec{y} \to \vec{0}} 0,
\end{multline*}
where the constant $C$ comes from the assumptions on the function $\phi$. Therefore, the measure $\mu$ would have to be the zero-measure, an impossibility given that the function $F$ is not identically zero.  

As a \HN function, the function $h$ can be written as
$$h(z_1,z_2) - \bar{h(w_1,w_2)} = \frac{1}{2\,\I}(z_1-\bar{w}_1)(z_2-\bar{w}_2)D(\vec{z},\vec{w}),$$
where the function $D$ is as specified in Theorem \ref{thm:positive_semi_decomposition} and equals
$$D((z_1,z_2),(w_1,w_2)) := \frac{-4}{(z_1-\bar{w}_1)(z_2-\bar{w}_2)}.$$
Let us now attempt to rewrite the above decomposition as
$$h(z_1,z_2) - \bar{h(w_1,w_2)} = (z_1-\bar{w}_1)\cdot\underbrace{\left(\frac{z_2-\bar{w}_2}{2\,\I}D(\vec{z},\vec{w})\right)}_{:= F_1(\vec{z},\vec{w})} + (z_2-\bar{w}_2)\cdot 0.$$
The function $F_1$, defined in the above way, can be shown to be equal to
$$F_1((z_1,z_2),(w_1,w_2)) = \frac{2\,\I}{z_1-\bar{w_1}},$$
providing, thusly, a valid decomposition of the form \eqref{eq:lowner_functions}.
\hfill$\lozenge$
\end{example}

\begin{example}
Let $h(z_1,z_2) := -(z_1+z_2)^{-1}$ for all $(z_1,z_2) \in \C^{+2}$. Then, it holds that
\begin{multline*}
    h(z_1,z_2) - \bar{h(w_1,w_2)} = -\frac{1}{z_1+z_2} + \frac{1}{\bar{w}_1 + \bar{w}_2} = \frac{z_1+z_2-\bar{w}_1 + \bar{w}_2}{(z_1+z_2)(\bar{w}_1+ \bar{w}_2)} \\
    = (z_1 - \bar{w}_1)\,\frac{1}{(z_1+z_2)(\bar{w}_1+ \bar{w}_2)} + (z_2 - \bar{w}_2)\,\frac{1}{(z_1+z_2)(\bar{w}_1+ \bar{w}_2)}
\end{multline*}
for all $\vec{z},\vec{w} \in \C^{+2}$. The function
$$F\colon ((z_1,z_2),(w_1,w_2)) \mapsto \frac{1}{(z_1+z_2)(\bar{w}_1+ \bar{w}_2)}$$
is indeed positive semi-definite on $\C^{+2} \times \C^{+2}$, as it is of the same form as the functions in Example \ref{ex:psd_simple}. However, it can be shown via an analogous reasoning as in the previous example that this functions is not a Poisson-type function.

As a \HN function, it admits a decomposition of the form
\begin{multline*}
h(z_1,z_2) - \bar{h(w_1,w_2)} \\ = \frac{1}{2\,\I}(z_1-\bar{w}_1)(z_2-\bar{w}_2)\cdot\underbrace{\frac{2\,\I\,(z_1+z_2-\bar{w}_1 -\bar{w}_2)}{(z_1-\bar{w}_1)(z_2-\bar{w}_2)(z_1+z_2)(\bar{w}_1+\bar{w}_2)}}_{= D(\vec{z},\vec{w})}
\end{multline*}
for all $\vec{z},\vec{w} \in \C^{+2}$. Investigating whether this decomposition can be rewritten into a decomposition of the form \eqref{eq:lowner_functions} as in the previous examples leads us to consider the function
$$(\vec{z},\vec{w}) \mapsto \frac{z_1 - \bar{w}_1}{2\,\I}\,D(\vec{z},\vec{w}) = \frac{z_1+z_2-\bar{w}_1 -\bar{w}_2}{(z_2-\bar{w}_2)(z_1+z_2)(\bar{w}_1+\bar{w}_2)}.$$
However, this function is no longer positive semi-definite. Indeed, choose $m = 2$, $\vec{z}_1 = (\I,-2+\I)$, $\vec{z}_2 = (2+\I,2\,\I)$, $c_2 = 1$ and
$$c_1 = \frac{\sqrt{377}-7}{2132}(-33-113\,\I).$$
Then, it holds that
$$\sum_{i,j=1}^2\frac{(\vec{z}_i)_1 - \bar{(\vec{z}_j)}_1}{2\,\I}D(\vec{z}_i,\vec{z}_j)c_i\bar{c}_j = \frac{4901-255\sqrt{377}}{8528} \cong -0.00588701 < 0.$$
Hence, a decomposition of the form \eqref{eq:lowner_functions} cannot always be constructed form a decomposition of the form \eqref{eq:positive_semi_decomposition} via the process that was exhibited in the previous example.
\hfill$\lozenge$
\end{example}

\section{Holomorphic functions on the unit polydisk with non-negative real part}
\label{sec:polydisk}

Historically, the idea to investigate the difference (or sum) of the values of a holomorphic function with a prescribed sign of its imaginary or real part dates back to the 1910's. A fundamental result from this era is due to Pick \cite[pg. 8]{Pick1915}, who proved, in modern terms, that for any holomorphic function $f\colon\D \to \C$ with non-negative real part it holds that
$$\frac{f(\xi) + \bar{f(\eta)}}{1-\xi\,\bar{\eta}} = \frac{1}{\pi}\int_{[0,2\pi)}\frac{1}{(\E^{\I\,s}-\xi)(\E^{-\I\,s}-\bar{\eta})}\diff\nu(s),$$
where $\xi,\eta \in \D$ and $\nu$ is the representing measure of the function $f$ in the sense of the Riesz-Herglotz representation theorem \cite[pg. 76]{ShilovGurevich1977}. Using the Cayley transform and its inverse, recalled explicitly later on, this result may be reinterpreted in the language of \HN functions where it reproduces the results summarized in the beginning of Section \ref{sec:main_thm}.

When considering functions of several variables, Kor\'anyi and Puk\'anszky gave an integral representation for holomorphic functions on the unit polydisk with non-negative real part \cite[Thm. 1]{KoranyiPukanszky1963}. Additionally, they gave a criterion to determine when a function defined \emph{a priori} on a subset of $\D^n$ of a particular type can be extended to holomorphic function on the whole of $\D^n$ with non-negative real part \cite[Thm. 2]{KoranyiPukanszky1963}. The latter condition of the latter theorem is expressed in terms of a certain function being positive semi-definite. Interestingly, this positive semi-definite function turns out to be a real constant multiple of the one appearing in the following adaptation of Theorem \ref{thm:positive_semi_decomposition}.

\begin{prop}
\label{prop:szego_kernel}
Let $f\colon \D^n \to \C$ be a holomorphic function. Then, $f$ has non-negative real part if and only if the 
\begin{equation}
    \label{eq:szego_kernel}
    (\vec{\xi},\vec{\eta}) \mapsto \frac{2^{n-1}}{\prod_{j=1}^n(1-\xi_j\bar{\eta}_j)}\cdot\big(f(\vec{\xi}) + \bar{f(\vec{\eta})}\big)
\end{equation}
is positive semi-definite on $\D^n \times \D^n$.
\end{prop}

\begin{remark}
The normalizing factor $2^{n-1}$ has no impact on the conclusion of the proposition, unlike the factor $(2\,\I)^{n-1}$ appearing in Theorem \ref{thm:Nevan_kernel_Nvar}. Furthermore, the function
$$(\vec{\xi},\vec{\eta}) \mapsto \prod_{j=1}^n\frac{1}{1-\xi_j\bar{\eta}_j},$$
where $\vec{\xi},\vec{\eta} \in \D^n$, is sometimes referred to as the \emph{Szeg\H{o} kernel} of $\D^n$, \cf \cite[pg. 450]{KoranyiPukanszky1963}.
\end{remark}

\proof
Recall first that the Cayley transform $\varphi\colon \C \to \C^+$ is defined as
$$\varphi(\zeta) = \I\,\frac{1+\zeta}{1-\zeta},$$
while its inverse $\varphi^{-1}\colon\C^+ \to \D$ is given by
$$\varphi^{-1}(z) = \frac{z - \I}{z + \I}.$$
Furthermore, if a change of variables between $s \in (0,2\pi)$ and $t \in \R$ is given by $\E^{\I\,s} = \varphi^{-1}(t)$, then
$$\diff s = \frac{2}{1+t^2}\diff t \quad\text{or}\quad \diff t = \frac{1}{1-\cos(s)}\diff s.$$

Let now $f\colon\D^n \to \C$ be a holomorphic function with non-negative real part. Then, the function
$$q(\vec{z}) := \I\,f(\varphi^{-1}(z_1),\ldots,\varphi^{-1}(z_n))$$ is a \HN function.
By Theorem \ref{thm:positive_semi_decomposition}, it holds that
\begin{multline}
\label{eq:temp2}
    q(\vec{z}) - \bar{q(\vec{w})} = \sum_{j=1}^nb_j(z_j - \bar{w}_j) \\ + \frac{1}{(2\,\I)^{n-1}}\prod_{j=1}^n(z_j - \bar{w}_j)\cdot\frac{1}{\pi^n}\int_{\R^n}\prod_{j=1}^n\frac{1}{(t_j - z_j)(t_j - \bar{w}_j)}\diff\mu(\vec{t}),
\end{multline}
where $\vec{b}$ and $\mu$ are the representing parameters of the function $q$ in the sense of Theorem \ref{thm:intRep_Nvar}.
Define now $\vec{\xi},\vec{\eta} \in \D^n$ and $\vec{s} \in (0,2\pi)^n$ via $\xi_j := \varphi^{-1}(z_j)$, $\eta_j := \varphi^{-1}(w_j)$ and $\E^{\I\,s_j} := \varphi^{-1}(t_j)$ for $j=1,2,\ldots,n$. In this notation, $f(\vec{\xi}) = -\I\,q(\vec{z})$ and equality \eqref{eq:temp2} becomes
\begin{multline}
    \label{eq:temp3}
    \I\,f(\vec{\xi}) + \I\,\bar{f(\vec{\eta})} \\
    = \I\,\sum_{j=1}^nb_j\left(\frac{1+\xi_j}{1-\xi_j} + \frac{1+\bar{\eta}_j}{1-\bar{\eta}_j}\right) + \frac{\I^n}{(2\,\I)^{n-1}}\prod_{j=1}^n\left(\frac{1+\xi_j}{1-\xi_j} + \frac{1+\bar{\eta}_j}{1-\bar{\eta}_j}\right) \\ 
    \cdot\frac{1}{(2\pi)^n}\int_{(0,2\pi)^n}\prod_{j=1}^n\frac{(1-\xi_j)(1-\bar{\eta}_j)}{(\E^{\I\,s_j}-\xi_j)(\E^{-\I\,s_j}-\bar{\eta}_j)}\diff\nu(\vec{s}) \\
    = 2\,\I\,\sum_{j=1}^nb_j\frac{1-\xi_j\bar{\eta}_j}{(1-\xi_j)(1-\bar{\eta}_j)} + \frac{\I}{2^{n-1}}\prod_{j=1}^n(1-\xi_j\bar{\eta}_j) \\
    \cdot\frac{1}{\pi^n}\int_{(0,2\pi)^n}\prod_{j=1}^n\frac{1}{(\E^{\I\,s_j}-\xi_j)(\E^{-\I\,s_j}-\bar{\eta}_j)}\diff\nu(\vec{s}).
\end{multline}
Here, the measure $\nu$ on $(0,2\pi)^n$ is a reparametrization of the measure $\mu$ obtained by setting
$$\diff\mu(\vec{t}) = \prod_{j=1}^n\frac{1}{1-\cos(s_j)}\diff\nu(\vec{s}).$$
Note now that by standard residue calculus it holds for every $\xi,\eta \in \D$ that
$$\frac{1}{\pi}\int_{[0,2\pi)}\frac{1}{(\E^{\I\,s}-\xi)(\E^{-\I\,s}-\bar{\eta})}\diff s = \frac{2}{1-\xi\,\bar{\eta}}.$$
Hence, for every $j =1,2\ldots,n$, we may write
\begin{multline*}
    2\,\I\,b_j\,\frac{1-\xi_j\bar{\eta}_j}{(1-\xi_j)(1-\bar{\eta}_j)} = \frac{\I}{2^{n-1}}\prod_{j=1}^n(1-\xi_j\bar{\eta}_j) \\ \cdot\frac{1}{\pi^n}\int_{[0,2\pi)^n}\prod_{j=1}^n\frac{1}{(\E^{\I\,s_j}-\xi_j)(\E^{-\I\,s_j}-\bar{\eta}_j)}\diff\sigma_j(\vec{s}),
\end{multline*}
where $\sigma_j$ is a measure on $[0,2\pi)^n$ defined as
$$\sigma_j := \lambda_{[0,2\pi)} \times \ldots \times \lambda_{[0,2\pi)} \times \underbrace{2\pi b_j \delta_0}_{j-\text{th coordinate}} \times \lambda_{[0,2\pi)} \times \ldots \times \lambda_{[0,2\pi)}.$$
Here, $\lambda_{[0,2\pi)}$ denotes the Lebesgue measure on ${[0,2\pi)}$ and $\delta_0$ denotes the Dirac measure at zero. We also extend the measure $\nu$ to a measure $\til{\nu}$ on $[0,2\pi)^n$ by setting
$$\til{\nu}|_{(0,2\pi)^n} := \nu \quad\text{and}\quad \til{\nu}|_{[0,2\pi)^n\setminus(0,2\pi)^n} := 0.$$
Denote now
$$\til{\sigma} := \til{\nu} + \sum_{j=1}^n\sigma_j.$$
Using this, equality \eqref{eq:temp3} may be rewritten as
$$f(\vec{\xi}) + \bar{f(\vec{\eta})} = \frac{1}{2^{n-1}}\prod_{j=1}^n(1-\xi_j\bar{\eta}_j)\cdot\frac{1}{\pi^n}\int_{[0,2\pi)^n}\prod_{j=1}^n\frac{1}{(\E^{\I\,s_j}-\xi_j)(\E^{-\I\,s_j}-\bar{\eta}_j)}\diff\til{\sigma}(\vec{s}).$$
Note finally that the function
\begin{equation}
    \label{eq:temp4}
    (\vec{\xi},\vec{\eta}) \mapsto \frac{1}{\pi^n}\int_{[0,2\pi)^n}\prod_{j=1}^n\frac{1}{(\E^{\I\,s_j}-\xi_j)(\E^{-\I\,s_j}-\bar{\eta}_j)}\diff\til{\sigma}(\vec{s})
\end{equation}
is positive semi-definite on $\D^n \times \D^n$, which follows by an analogous calculation as was done in the proof of Lemma \ref{lem:positive_semi_poisson}. This finishes the first part of the proof.

Conversely, assume that we have a function $f\colon \D^n \to \C$ for which the function \eqref{eq:szego_kernel} is positive semi-definite on $\D^n \times \D^n$. To show that $f$ has non-negative real part, we only need to evaluate the function \eqref{eq:szego_kernel} at $\vec{\xi} = \vec{\eta}$ and use Lemma \ref{lem:positive_semi_def_properties}. This finishes the proof.
\endproof

\begin{remark}
The Nevanlinna kernel and the function \eqref{eq:szego_kernel} are not equivalent under the Cayley transform. More precisely, let $\varphi$ and $\varphi^{-1}$ be the Cayley transform and its inverse as in the proof of Proposition \ref{prop:szego_kernel}, let $\vec{z},\vec{w} \in \C^{+n}$ and $\vec{\xi},\vec{\eta} \in \D^n$ be such that $\xi_j := \varphi^{-1}(z_j)$ and $\eta_j := \varphi^{-1}(w_j)$ and let $q\colon\C^{+n} \to \C$ be a \HN function and $f\colon\D^n \to \C$ be a holomorphic function with non-negative real part such that $q(\vec{z}) = \I\,f(\vec{\xi})$. Then,
$$\frac{(2\,\I)^{n-1}}{\prod_{j=1}^n(z_j - \bar{w}_j)}\big(q(\vec{z}) - \bar{q(\vec{w})}\big) \neq \frac{2^{n-1}}{\prod_{j=1}^n(1-\xi_j\bar{\eta}_j)}\cdot\big(f(\vec{\xi}) + \bar{f(\vec{\eta})}\big).$$
An analogous non-equivalence under the Cayley transform holds for Poisson-type function and function of the type \eqref{eq:temp4}.
\end{remark}

\section*{Acknowledgements}

The author would like to thank Dale Frymark for many enthusiastic discussions on the subject.

\bibliographystyle{amsplain}
\bibliography{MitjaNedic_characterization_via_psd_functions}

\end{document}